\documentclass[a4paper,11pt]{article}

\usepackage{latexsym,amsfonts,amssymb,amsmath,amscd}
\usepackage[utf8]{inputenc}

\newtheorem{theorem}{Theorem}[section]
\newtheorem{lemma}[theorem]{Lemma}

\newtheorem{corollary}[theorem]{Corollary}

\newtheorem{example}[theorem]{Example}

\newtheorem{remark}[theorem]{Remark}

\date{}

\begin{document}



\parindent    = 25pt
\baselineskip = 16pt

\mbox{}\\
\begin{center}
{\Large\bf Genus expanded cut-and-join operators and generalized Hurwtiz numbers}
\vskip 0.4in

Quan Zheng\footnote{E-mail: quanzheng2100@163.com. Partially supported by NSFC 11171258 and NSFC 11571242.}\qquad

 Mathematics College, Sichuan University\\
610064, Chengdu, Sichuan, PRC

November 1, 2015\\
(Dedicated to Prof. An-Min Li for his 70th birthday)
\vskip 0.4in
{\bf Abstract}\\[\baselineskip]

\parbox{14cm}{\quad
  To distinguish the contributions to the generalized Hurwitz number  of the source Riemann surface with different genus, by observing carefully the symplectic surgery and the gluing formulas of the relative GW-invariants, we define the genus expanded  cut-and-join operators. Moreover all normalized the genus expanded  cut-and-join operators with same degree form  a differential algebra, which is isomorphic to the central subalgebra of the symmetric group algebra.
  As an application, we get some  differential equations for the generating functions of the generalized
  Hurwitz numbers  for the source Riemann surface with different genus,
  thus we can express the generating functions in  terms of the genus expanded  cut-and-join operators.
}
\end{center}
{\bf Keyword:} genus expanded cut-and-join operator; differential algebra; Hurwitz number; generating function\\
{\bf Subject Classification: 14N10, 14N35, 05E10}
\indexspace


\renewcommand\contentsname{\Large Contents}

\tableofcontents

\section{Introduction}

 Hurwitz Enumeration Problem has been studied for more than one century,
 which has been extensively  applied to mathematics and physics.
 It is an interesting project for mathematic physics,  geometry, representation theory,
 integrable  systems, Hodge integral,  etc.
 Hurwitz numbers have many different expressions in the different fields, for example,
 \cite{[AMMN]}, \cite{[DLM]}, \cite{[ELSV]}, \cite{[GJ1]}, \cite{[GJ2]},\cite{[GJ3]}, \cite{[H]}, \cite{[MMN]}, \cite{[MMN1]}, \cite{[MMN2]}.
 One of the important geometric tools to deal with Hurwitz numbers is the so-called symplectic surgery:
 cutting and gluing \cite{[IP]}, \cite{[L]},\cite{[LR]}, in the views of algebra and differential equations,
 which is equivalent to the so-called cut-and-join operators \cite{[GJ1]}, \cite{[GJ2]}, \cite{[GJ3]},\cite{[LZZ]}, \cite{[OP]}.
 The standard cut-and-join operators can be used to deal with the almost simple Hurwitz numbers
 and the almost simple double Hurwitz numbers  \cite{[GJ1]}, \cite{[GJ2]},\cite{[GJ3]}, \cite{[LZZ]}, \cite{[OP]}.
 A. Mironov, A. Morozov, and S. Natanzon, etc, have defined the generalized
 cut-and-join operators in terms of the matrix Miwa variable \cite{[AMMN]}, \cite{[MMN]}, \cite{[MMN1]}.
 However, the generating functions obtained by the generalized cut-and-join operators
 defined by A. Mironov, A. Morozov, and S. Natanzon cannot distinguish the contributions of
 the source Riemann surface with the different genus, which is obviously  interesting and important.
 In our paper,  by observing carefully the symplectic surgery and the gluing formulas of the
 relative GW-invariants developed by A.M. Li and Y.B. Ruan \cite{[LR]},  and by E. Ionel  and T. Parker  \cite{[IP]},
we define the genus expanded  cut-and-join operators to distinguish the contributions of
the source Riemann surfaces with different genus. A little more precisely, we introduce one parameter $z$
to ``mark" the genus,  then   we associate to
 every partition $\Delta$   a genus  expanded cut-and-join operator $W(\Delta,z)$ as follows:
\begin{equation}\label{4}
W(\Delta,z)=\sum_{\Gamma,\Gamma^{\prime}\vdash d}
\frac{d!}{|C_{\Gamma^{\prime}}|} z^{d+l(\Gamma^{\prime})-l(\Delta)
-l(\Gamma)}\mu_{0}^{h_+,d}(\Gamma^{\prime},\Delta,\Gamma)
p_{\Gamma}\frac{\partial}{\partial p_{\Gamma^{\prime}}}.
\end{equation}
For the  precise meaning of the notations in the above formula, see the section 2 and the section 3.
 Moreover, the power of $z$ is exactly the ``lost" genus after we have
executed the symplectic cutting on the target Riemann surface, referring to  Remark 4.3. \\

One of our main theorems is the followings, referring the theorem 3.5:\\
\bigskip
{\bf Theorem A}
For a given $d$, as operators on the functions of the time-variables
$p = (p_1, p_2,\cdots )$, we have
\begin{equation}
 {W(\Delta_1,z) W(\Delta_2,z) = \sum_{\Delta_3}z^{d-l(\Delta_1)-l(\Delta_2)+l(\Delta_3)}
C_{\Delta_1\Delta_2}^{\Delta_3}  W(\Delta_3,z) },
\end{equation}
where
\begin{equation}
C_{\Delta_1\Delta_2}^{\Delta_3}=\frac{d!}{|C_{\Delta_3}|}\mu_{0}^{h,d}(\Delta_{1},\Delta_{2},\Delta_{3})
\end{equation}
is the  structure constants
of the center subalgebra of the group algebra $\mathbb{C}[S_{d}]$.

A more interested result in this paper is that all normalized the genus expanded  cut-and-join operators with same degree form  a differential algebra, which is isomorphic to the central subalgebra of the symmetric group algebra, referring Corollary (\ref{cor 3.6}).\\

As an application of the genus expanded cut-and-join operators, we define the
generating function as (referring the section 4)
\begin{eqnarray}
&&\Phi_g\{z|(u_1,\Delta_1),\cdots,(u_n,\Delta_n)|p^{(1)},\cdots,p^{(k)},p\} \nonumber\\
&=& \sum_{l_1,\cdots,l_n\ge 0}
\sum_{\Gamma,\Gamma_1,\ldots,\Gamma_{k}}z^{2h-2}\mu_{g}^{h,d}(\underbrace{\Delta_{1},\cdots,\Delta_1}_{l_1},\cdots,\underbrace{\Delta_{n},\cdots,\Delta_n}_{l_n},\Gamma_1,\cdots,
\Gamma_k,\Gamma)[\prod_{j=1}^{n}\frac{u_j^{l_j}}{l_j!}][\prod_{i=1}^{k}p^{(i)}_{\Gamma_i}]p_\Gamma \nonumber\\
&=&\sum_{l_1,\cdots,l_n\ge 0}\sum_{\Gamma,\Gamma_1,\ldots,\Gamma_{k}}z^{2h-2}({\frac {dim{\lambda}}{d!}})^{2-2g}[\prod_{\makebox {j=1}}^n(\phi_{\lambda}(\Delta_j))^{l_j}\frac{(u_j)^{l_j}}{l_j!}]
\phi_{\lambda}(\Gamma)[\prod_{i=1}^{k}\phi_{\lambda}(\Gamma_i)p^{(i)}_{\Gamma_i}]p_\Gamma.
\end{eqnarray}
Then we have our another main theorem (referring theorem 4.2):\\
{\bf Theorem B} For any $i$, we have
\begin{eqnarray}
&&\nonumber \frac{\partial\Phi_g\{z|(u_1,\Delta_1),\cdots,(u_n,\Delta_n)|p^{(1)},\cdots,p^{(k)},p\}}{\partial u_i} \nonumber\\
& =&W(\Delta_i,z)\Phi_g\{z|(u_1,\Delta_1),\cdots,(u_n,\Delta_n)|p^{(1)},\cdots,p^{(k)},p\}.
\end{eqnarray}
Thus we can express the generating functions in  terms of the genus expanded  cut-and-join operators
once we know the initial values. Moreover,  the generating functions can be directly expressed
by the  expansion in terms of the genus of the source Riemann surface, referring to Example 4.5.

Our paper is deeply inspired by A. Mironov, A. Morozov, and S. Natanzon's works \cite{[AMMN]}, \cite{[MMN]}, \cite{[MMN1]}, \cite{[MMN2]}.

\section{Generalized Hurwitz numbers}
Let $\Sigma^h$ be a compact (possibly  disconnected) Riemann surface of genus $h,$ and $\Sigma^g$
 a compact connected Riemann surface of genus $g$.
 For a given point set $\{q_1,\cdots,q_k\}\in \Sigma^g$, which is called the set of
 branch points, we call a
holomorphic map
$f:\Sigma^h\to \Sigma^g$
a ramified
covering of $\Sigma^g$ of degree $d\ge 0$ by $\Sigma^h$ with a ramification type $(\Delta_{1},\cdots,\Delta_{k})$, if
the preimages of $f^{-1}(q_i)=\{p_i^{1},\cdots,p_i^{m_i}\}$  with orders $\Delta_i = (\delta^i_1,\cdots,\delta^i_{m_i})$ for $i=1,\cdots,k,$ respectively.
Two ramified coverings $f_1$  and $f_2$ with type $(\Delta_{1},\cdots,\Delta_{k})$
 are said to be equivalent if there is a
homeomorphism
$\pi: \Sigma^h\to \Sigma^h$ such that $f_1=f_2\circ\pi$ and $\pi$
preserves the preimages and the ramification type of $f_1$ and $f_2$ at each point $q_i\in \Sigma^g$.
Let $\mu_{g}^{h,d}(\Delta_{1},\cdots,\Delta_{k})$ be the number of equivalent
 covering of $\Sigma^g$ by $\Sigma^h$ with ramification type $(\Delta_{1},\cdots,\Delta_{k})$,
 which is called the generalized Hurwitz number.
  Note that $\mu_{g}^{h,d}(\Delta_{1},\cdots,\Delta_{k})$ is nonzero only if the Hurwitz formula
\begin{equation}\label{10}
(2-2g)d-(2-2h)=\sum_{i=1}^{k}(d-m_i)
\end{equation}
holds,
where $m_i$ is the numbers of the preimage of $q_i$. How to determine $\mu_{g}^{h,d}(\Delta_{1},\cdots,\Delta_{k})$
is known as the generalized
Hurwitz Enumeration Problem \cite{[H]}.
In general, $\Delta_i$ is not required  to be a partition of $d$, for simplicity, throughout  this paper,
we will assume  it. Denote by $\Delta_i \vdash d$ if $\Delta_i$ is a partition of $d.$

\bigskip
Let  $S_d$ be the symmetric group of $d$ letters, $\mathbb{C}[S_d]$ be the group algebra of $S_d$ and $Z\mathbb{C}[S_d]$ be the central subalgebra of $\mathbb{C}[S_d]$,
then in the terms of algebra, the story begins from a simple fact as the following formula \cite{[H]},
\begin{eqnarray}
\mu_{g}^{h,d}(\Delta_{1},\cdots,\Delta_{k})=\frac{1}{d!}[1]
\prod_{j=1}^g\prod_{a_j,b_j\in S_d}[a_j,b_j]C_{\Delta_1}\cdots C_{\Delta_k},
\end{eqnarray}
where $[1]$ means that we take the coefficient of the identity  of the product of $g$-tuple  commutators $\prod_{a,b\in S_d}[a,b]$ and $k$-tuple
  central elements $C_{\Delta_1}, \cdots, C_{\Delta_k}\in Z\mathbb{C}[S_d] $
  corresponding to the partitions $\Delta_1, \cdots, \Delta_k.$
  Thus, according to the theory of representations of the symmetric group $S_d$, we arrive at the formula \cite{[OP]},
\begin{eqnarray}{}
\mu_{g}^{h,d}(\Delta_{1},\cdots,\Delta_{k})=
\sum_{\lambda \vdash d}({\frac {dim{\lambda}}{d!}})^{2-2g}\phi_{\lambda}(\Delta_1)\cdots\phi_{\lambda}(\Delta_k),
\end{eqnarray}
which expresses them through the properly normalized symmetric group characters
\begin{equation}\label{21}
\phi_{\lambda}(\Delta)=\frac{1}{dim\lambda}|C_\Delta|\chi_\lambda(\Delta),
\end{equation}
where   $\lambda$ is a Young diagram of degree $d$, and $dim\lambda$ is the dimension of the irreducible  representation of the symmetric group $S_{d}$ corresponding to $\lambda$, $\chi_\lambda(\Delta)$ is the  character of a permutation $\sigma\in C_{\Delta}$ under the irreducible  representation $\lambda$, moreover, $|C_\Delta|$ is  the number of the permutations of $S_d$ with cyclic type $\Delta$.

Let $p=(p_1,p_2,p_3,\cdots,)$ be indeterminantes, which are called the
time-variables, and we assume that $\Delta=(\delta_{1},\cdots,\delta_n)$ is a partition. We denote
 $$l(\Delta):=n,$$
$$m_r(\Delta):=\#\{i|\delta_i=r,\delta_i \in \Delta \},$$
$$\Delta!:=\prod_{r\ge 1}m_r(\Delta)!,$$
$$p_\Delta:=p_{\delta_1}\cdots p_{\delta_n},$$
and
$$\frac{\partial}{\partial p_\Delta}:=\frac{1}{\Delta!}\frac{\partial}{\partial p_{\delta_1}}\cdots \frac{\partial}{ \partial p_{\delta_n}}.$$

For any partition $\lambda \vdash d,$ it is well known that the symmetric group
characters $\phi_{\lambda}(\Delta)$ are related to the Schur functions $S_{\lambda}\{p\}$
as following \cite{[MMN]},
\begin{equation}\label{5}
S_{\lambda}\{p\} = \sum_{\Gamma^{\prime}}
\frac{ dim{\lambda}}{d!}\phi_{\lambda}(\Gamma^{\prime})p_{\Gamma^{\prime}}.
\end{equation}

 Our first important  lemma is about the associativity of the generalized Hurwitz number:

\lemma \label{lemma1.1} \cite{[MMN2]} For any positive integer $k> l\ge 1$, and $g=g_1+g_2,$ we have
\begin{equation}\label{6}
\mu_{g}^{h,d}(\Delta_{1},\cdots,\Delta_{k})=
 \sum_{\Delta^{\prime}}\mu_{g_1}^{h_1,d}(\Delta_{1},\cdots,\Delta_{l},\Delta^{\prime})\frac{d!}{|C_{\Delta^{\prime}}|}\mu_{g_2}^{h_2,d}(\Delta^{\prime},\Delta_{l+1},\cdots,\Delta_{k}).
\end{equation}

{\bf{Proof}}   Algebraic proof: by formula(3), we have
\begin{eqnarray*}
\mathrm{RHS}&=&\sum_{\Delta^{\prime}}\sum_{\lambda}
({\frac {dim{\lambda}}{d!}})^{2-2g_1}\phi_{\lambda}(\Delta_1)\cdots\phi_{\lambda}(\Delta_l)
\phi_{\lambda}(\Delta^{\prime})\frac{d!}{|C_{\Delta^{\prime}}|}\\
&&\times
\sum_{\kappa}({\frac {dim{\kappa}}{d!}})^{2-2g_2}\phi_{\kappa}(\Delta^{\prime})
\phi_{\kappa}(\Delta_{l+1})\cdots\phi_{\kappa}(\Delta_k).
\end{eqnarray*}
Then the lemma follows from the orthogonal  relation of the irreducible characteristic  of symmetric group  $S_d$:
\begin{equation}
\frac{1}{d!}\sum_{\Delta^{\prime}}\chi_{\lambda}(\Delta^{\prime})|C_{\Delta^{\prime}}|
\chi_{\kappa}(\Delta^{\prime})=\delta_{\lambda,\kappa}.
\end{equation}

Geometric proof: it is well known \cite{[LZZ]} that we can interpret the generalized
Hurwitz numbers $\mu_{g}^{h,d}(\Delta_{1},\cdots,\Delta_{k})$ as the  relative
Gromov-Witten invariants defined by  A.M. Li and Y.B. Ruan \cite{[LR]}, and by E. Ionel  and T. Parker \cite{[IP]},
thus we only need to execute symplectic surgery \cite{[L]}, \cite{[LR]} by cutting the Riemannian surface $\Sigma_g$
into two parts:  one part has the branch point $P_1, \cdots,P_l$ with genus $g_1$, and another part has
the rest branch points $P_{l+1}, \cdots,P_k$ with genus $g_2$, moreover each part has one more branch
point $P^{\prime}$, which corresponds to  the infinity end \cite{[LR]}, with the ramification type (or  tangent multiplicities \cite{[LR]})
$\Delta^{\prime}$. Note that $g=g_1+g_2.$
Thus we have to consider two relative GW-invariants:
$\mu_{g_1}^{h_1,d}(\Delta_{1},\cdots,\Delta_{l},\Delta^{\prime})$ and
$\mu_{g_2}^{h_2,d}(\Delta^{\prime},\Delta_{l+1},\cdots,\Delta_{k})$,  which are nonzero only if
\begin{equation}\label{8}
(2-2g_1)d-(2-2h_1)=\sum_{i=1}^{l}(d-m_i)+(d-l(\Delta^{\prime})),
\end{equation}
or
\begin{equation}\label{9}
(2-2g_2)d-(2-2h_2)=(d-l(\Delta^{\prime}))+\sum_{i=l+1}^{k}(d-m_i),
\end{equation}
respectively.
The factor $\frac{d!}{|C_{\Delta^{\prime}}|}$ in formula (\ref{6}) comes from the tangent multiplicities
and the automorphisms at the infinity end \cite{[LR]}.
Then applying
the  gluing formula \cite{[LR]}, mimicking what we have done in \cite{[LZZ]}, we can  derive the above lemma.

\section{Genus expanded  cut-and-join operators}
Observing the symplectic surgery and the gluing formula [LR], we  can associate to
 every partition $\Delta \vdash d$   a genus  expanded cut-and-join operator $W(\Delta,z)$ as follows:
\begin{equation}\label{4}
W(\Delta,z)=\sum_{\Gamma,\Gamma^{\prime}\vdash d}
\frac{d!}{|C_{\Gamma^{\prime}}|} z^{d+l(\Gamma^{\prime})-l(\Delta)
-l(\Gamma)}\mu_{0}^{h_+,d}(\Gamma^{\prime},\Delta,\Gamma)
p_{\Gamma}\frac{\partial}{\partial p_{\Gamma^{\prime}}},
\end{equation}
where  genus $h_+$ is determined  by the Hurwitz formula:
\begin{equation}
2d-(2-2h_+)=(d-l(\Gamma^{\prime}))+(d-l(\Delta))+(d-l(\Gamma)).
\end{equation}
\bigskip
\remark
The power $d+l(\Gamma^{\prime})-l(\Delta)-l(\Gamma)$ of $z$ in the formula (\ref{4})
is the ``lost" genus after we have executed  the symplectic surgery on the Riemann surface,
referring  to Remark 4.3, which is the reason for us to define the genus expanded cut-and-join operator like
formula (\ref{4}), rather than the normalized genus expanded cut-and-join operator as formula (\ref{20}).

\bigskip
\remark  A. Mironov, A. Morozov, and S. Natanzon, etc,
have defined the cut-and-join operators in terms of the matrix Miwa variable \cite{[AMMN]},
\cite{[MMN]}, \cite{[MMN1]}, but it is easy to check that these two definitions
coincide if we take $z=1$ and keep the degree $d$ invariant.

\bigskip
To find out the eigenfunctions of the genus expanded cut-and-join operator  $W(\Delta,z)$,
we define the ``genus expanded"
Schur functions $S_{\lambda}\{p,z\}(\lambda \vdash d)$ similar to formula (\ref{5})
as follows:
\begin{equation}\label{11}
S_{\lambda}\{p,z\} := \sum_{\Gamma^{\prime}} z^{-d-l(\Gamma^{\prime})}\frac{ dim{\lambda}}{d!}\phi_{\lambda}(\Gamma^{\prime})p_{\Gamma^{\prime}}
\end{equation}

First of all,
 for any $\Delta,\Delta^{\prime} \vdash d,$ we note that
\begin{equation}
\frac{\partial}{\partial p_\Delta}p_{\Delta^{\prime}}=\delta_{\Delta,\Delta^{\prime}}.
\end{equation}
  Then by formulas (\ref{4}), (\ref{11}),
 we have
\lemma For any partition $\Gamma^{\prime}\vdash d,$ the following equality holds:
\begin{equation}\label{14}
W(\Delta,z)p_{\Gamma^{\prime}}=\sum_{\Gamma\vdash d}\frac{d!}{|C_{\Gamma^{\prime}}|}z^{d+l(\Gamma^{\prime})-l(\Delta)-l(\Gamma)}\mu_{0}^{h_+,d}(\Gamma^{\prime},\Delta,\Gamma)p_{\Gamma}. \end{equation}
Moreover $W(\Delta,z)$ have the genus expanded Schur function  $S_\lambda\{p,z\}$
 as their eigenfunctions
and $z^{d-l(\Delta)}\phi_\lambda(\Delta)$ as the corresponding eigenvalues:
\begin{equation}\label{12}
 W(\Delta,z)S_\lambda\{p,z\} = z^{d-l(\Delta)}\phi_\lambda(\Delta) S_\lambda\{p,z\}.
\end{equation}
?Proof By straightforward  calculation, we omit it.

\bigskip
\example
\begin{itemize}
\item[(1)] In the simplest
case $\Delta=(2)\vdash 2,$ we have
$$
W(\Delta,z)=\frac{1}{2}z^2p_2\frac{\partial^2}{\partial p_1\partial p_1}+p_1p_1\frac{\partial}{\partial p_2},
$$
which is the standard cut-and-join operator \cite{[GJ1]},\cite{[GJ2]},\cite{[GJ3]}.

\item[(2)]
For $\Delta=(2,1)\vdash 3,$ we have
$$
W(\Delta,z)=z^2[2p_3\frac{\partial^2}{\partial p_1\partial p_2}+\frac{1}{2}p_1p_2\frac{\partial^3}{\partial p_1\partial p_1\partial p_1}]+[p_1^3\frac{\partial^2}{\partial p_1\partial p_2}+3p_1p_2\frac{\partial}{\partial p_3}].
$$
\end{itemize}

\bigskip
\theorem
For a given $d$, as operators on the functions of the time-variables
$p = (p_1, p_2,\cdots )$, we have
\begin{equation}\label{19}
 {W(\Delta_1,z) W(\Delta_2,z) = \sum_{\Delta_3}z^{d-l(\Delta_1)-l(\Delta_2)+l(\Delta_3)}
C_{\Delta_1\Delta_2}^{\Delta_3}  W(\Delta_3,z) },
\end{equation}\label{CWW0}
where
\begin{equation}
C_{\Delta_1\Delta_2}^{\Delta_3}=\frac{d!}{|C_{\Delta_3}|}\mu_{0}^{h,d}(\Delta_{1},\Delta_{2},\Delta_{3})
\end{equation}
is the  structure constants
of the center subalgebra of the group algebra $\mathbb{C}[S_{d}]$.
?Proof
For any $\Gamma^{\prime}\vdash d$, by formula (\ref{14}), we have
\begin{eqnarray}\label{13}
&&W(\Delta_1,z) W(\Delta_2,z)p_{\Gamma^{\prime}}\nonumber\\
&=&W(\Delta_1,z)\sum_{\Gamma}\frac{d!}{|C_{\Gamma^{\prime}}|}z^{d+l(\Gamma^{\prime})-l(\Delta_2)-l(\Gamma)}\mu_{0}^{h_2,d}(\Gamma^{\prime},\Delta_2,\Gamma)p_{\Gamma}\nonumber\\
&=&\sum_{\Delta}z^{2d+l(\Gamma^{\prime})-l(\Delta_1)-l(\Delta_2)-l(\Delta)}\frac{d!}{|C_{\Gamma^{\prime}}|}\sum_{\Delta^{\prime}}\mu_{0}^{h_2,d}(\Gamma^{\prime},\Delta_2,\Delta^{\prime})\frac{d!}{|C_{\Delta^{\prime}}|}\mu_{0}^{h_1,d}(\Delta^{\prime},\Delta_1,\Delta)p_{\Delta} \nonumber\\
&=&\sum_{\Delta}z^{2d+l(\Gamma^{\prime})-l(\Delta_1)-l(\Delta_2)-l(\Delta)}\frac{d!}{|C_{\Gamma^{\prime}}|}\mu_{0}^{h^{\prime},d}(\Gamma^{\prime},\Delta_2,\Delta_1,\Delta)p_{\Delta}\nonumber\\
&=&\sum_{\Delta}z^{2d+l(\Gamma^{\prime})-l(\Delta_1)-l(\Delta_2)-l(\Delta)}\frac{d!}{|C_{\Gamma^{\prime}}|}\sum_{\Delta_3}\mu_{0}^{h,d}(\Delta_1,\Delta_2,\Delta_3)\frac{d!}{|C_{\Delta_3}|}\mu_{0}^{h_3,d}(\Delta_3,\Gamma^{\prime},\Delta)p_{\Delta}\nonumber\\
&=&\sum_{\Delta}z^{2d+l(\Gamma^{\prime})-l(\Delta_1)-l(\Delta_2)-l(\Delta)}\frac{d!}{|C_{\Gamma^{\prime}}|}\sum_{\Delta_3}C_{\Delta_1\Delta_2}^{\Delta_3}\mu_{0}^{h_3,d}(\Delta_3,\Gamma^{\prime},\Delta)p_{\Delta}\nonumber\\
&=&\sum_{\Delta_3}z^{d-l(\Delta_1)-l(\Delta_2)+l(\Delta_3)}C_{\Delta_1\Delta_2}^{\Delta_3}W(\Delta_3,z)p_{\Gamma^{\prime}},
\end{eqnarray}
which is equivalent to formula (\ref{19}), and where $h_1, h_2, h_3, h^{\prime},h$ are the genus.
From geometric viewpoint, the  third and fifth equalities
imply that we can obtain the same results  if we execute the symplectic surgery once instead of twice;
or from algebraic
viewpoint, the third and fifth equality are equivalent to  the associativity  of the generalized
Hurwitz numbers as in Lemma 1.1.
\bigskip
\corollary \label{cor 3.6}
If we normalize the genus expanded cut-and-join operator $W(\Delta, z)$ by a factor  $z^{-d+l(\Delta)}$
\begin{equation}\label{20}
\hat{W}(\Delta, z):=z^{-d+l(\Delta)}W(\Delta, z),
\end{equation}
then for a given $d$, as operators on the space of functions in time-variables $p=(p_1, p_2, \cdots)$,
all genus expanded  cut-and-join operators $\hat{W}(\Delta, z)$ for $\Delta \vdash d$ form a commutative associative
algebra, denoted by ${\cal{W}}_d,$
 \begin{equation}
 {\hat{W}(\Delta_1,z) \hat{W}(\Delta_2,z) = \sum_{\Delta_3}
C_{\Delta_1\Delta_2}^{\Delta_3} \hat{W}(\Delta_3,z) }
\end{equation}\label{CWW0}
 i.e., we have an algebraic isomorphism:
\begin{eqnarray}
{\cal{W}}_d &\cong &Z(\mathbb{C}[S_{d}])\nonumber\\
        \hat{W}(\Delta,z)   & \mapsto & C_{\Delta}.
\end{eqnarray}
Moreover, by formula (\ref{12}), $\hat{W}(\Delta,z)$ have the genus expanded Schur function  $S_\lambda\{p,z\}$
 as their eigenfunctions
and $\phi_\lambda(\Delta)$ as the corresponding eigenvalues:
\begin{equation}
 \hat{W}(\Delta,z)S_\lambda\{p,z\} = \phi_\lambda(\Delta) S_\lambda\{p,z\}.
\end{equation}
?Proof By straightforward  calculation, we omit it.

\section {Generating functions and its differential equations}
 For any given genus $g\ge 0,$ any  given degree $d$ and partitions $\Delta_1, \cdots, \Delta_n \vdash d,$ to apply the genus expanded cut-and-join operators (\ref{4}), we define the
generating function as
\begin{eqnarray}
&&\Phi_g\{z|(u_1,\Delta_1),\cdots,(u_n,\Delta_n)|p^{(1)},\cdots,p^{(k)},p\} \nonumber\\
&=& \sum_{l_1,\cdots,l_n\ge 0}
\sum_{\Gamma,\Gamma_1,\ldots,\Gamma_{k}}z^{2h-2}\mu_{g}^{h,d}(\underbrace{\Delta_{1},\cdots,\Delta_1}_{l_1},\cdots,\underbrace{\Delta_{n},\cdots,\Delta_n}_{l_n},\Gamma_1,\cdots,
\Gamma_k,\Gamma)[\prod_{j=1}^{n}\frac{u_j^{l_j}}{l_j!}][\prod_{i=1}^{k}p^{(i)}_{\Gamma_i}]p_\Gamma \nonumber\\
&=&\sum_{l_1,\cdots,l_n\ge 0}\sum_{\Gamma,\Gamma_1,\ldots,\Gamma_{k}}z^{2h-2}({\frac {dim{\lambda}}{d!}})^{2-2g}[\prod_{\makebox {j=1}}^n(\phi_{\lambda}(\Delta_j))^{l_j}\frac{(u_j)^{l_j}}{l_j!}]
\phi_{\lambda}(\Gamma)[\prod_{i=1}^{k}\phi_{\lambda}(\Gamma_i)p^{(i)}_{\Gamma_i}]p_\Gamma,
\end{eqnarray}
where $z, u_1,\cdots, u_n$ are indeterminate variables, $p, p^{(1)}, \cdots, p^{(k)}$ are time-variables, and
$2h-2$ is determined  by the Hurwitz formula:
\begin{equation}\label{g1}
(2-2g)d-(2-2h)=\sum_{j=1}^{n}l_j(d-l(\Delta_j))+\sum_{j=1}^{k}(d-l(\Gamma_j))+(d-l(\Gamma)).
\end{equation}
Moreover, we have some special initial values:
\item
\begin{eqnarray}\label{1}
\Phi_0\{z||p \}&=&\sum_{\lambda}\sum_{\Delta}z^{-d-l(\Delta)}(\frac{dim{\lambda}}{d!})^2\phi_{\lambda}(\Delta)p_{\Delta}\nonumber\\
&=&\sum_{\lambda}\frac{dim{\lambda}}{d!}S_{\lambda}\{p,z\}\nonumber\\
&=&z^{-2d}\frac{p_1^d}{d!};
\end{eqnarray}
\begin{eqnarray}\label{2}
\Phi_0\{z||p^{(1)},p\}&=&\sum_{\lambda}\sum_{\Delta_1,\Delta_2}z^{-l(\Delta_1)-l(\Delta_2)}(\frac{dim{\lambda}}{d!})^2\phi_{\lambda}(\Delta_1)\phi_{\lambda}(\Delta_2)p^{(1)}_{\Delta_1}p_{\Delta_2}\nonumber\\
&=&\sum_{\lambda}z^{2d}S_{\lambda}\{p^{(1)},z\}S_{\lambda}\{p,z\}\nonumber\\
&=&\sum_{\Delta}z^{-2l(\Delta)}\frac{|C_\Delta|}{d!}p^{(1)}_\Delta p_\Delta
\end{eqnarray}

\remark The formula (\ref{1}) can be regarded  as the genus and degree  expansion of hook formula \cite{[M]},\cite{[MMN]} and the formula (\ref{2}) as  the genus and degree expansion  of the Cauchy-Littlewood identity \cite{[M]},\cite{[MMN]}. Moreover one side of both of the formulas is algebraic expression of Hurwitz numbers and the another side is geometric expression of the Hurwitz numbers.

\bigskip
\theorem For any $i$, we have
\begin{eqnarray}\label{3}
&&\nonumber \frac{\partial\Phi_g\{z|(u_1,\Delta_1),\cdots,(u_n,\Delta_n)|p^{(1)},\cdots,p^{(k)},p\}}{\partial u_i} \nonumber\\
& =&W(\Delta_i,z)\Phi_g\{z|(u_1,\Delta_1),\cdots,(u_n,\Delta_n)|p^{(1)},\cdots,p^{(k)},p\}.
\end{eqnarray}
?Proof
Obviously, we have
$$
\frac{\partial\Phi_g\{z|(u_1,\Delta_1),\cdots,(u_n,\Delta_n)|p^{(1)},\cdots,p^{(k)},p\}}{\partial u_i}$$
$$=\sum_{l_1,\cdots,l_n\ge 0}
\sum_{\Gamma_1,\ldots,\Gamma_{k}, \Gamma}z^{2h-2}\mu_{g}^{h,d}(\underbrace{\Delta_{1},\cdots,\Delta_1}_{l_1},\cdots,\underbrace{\Delta_{n},\cdots,\Delta_n}_{l_n},\Gamma_1,\cdots,
\Gamma_k,\Gamma)$$
$$\frac{(u_i)^{l_i-1}}{(l_i-1)!}[\prod_{\makebox {j=1,}\\ j\ne i}^n\frac{(u_j)^{l_j}}{l_j!}]
[\prod_{j=1}^{k}p^{(j)}_{\Gamma_j}]p_\Gamma.
$$
We can write RHS of equation (\ref{3}) as
\begin{eqnarray*}
\mathrm{RHS}&=&\sum_{l_i\ge 1} \sum_{l_1,\cdots,\check{l_i},\cdots,l_n\ge 0}
\sum_{\Gamma_1,\ldots,\Gamma_{k},\Gamma^{\prime}}z^{2h^i_--2}\\
&&\times
\mu_{g}^{h^i_-,d}(\underbrace{\Delta_{1},\cdots,\Delta_1}_{l_1},\cdots,\underbrace{\Delta_{i},
\cdots,\Delta_i}_{l_i-1},\cdots,\underbrace{\Delta_{n},\cdots,\Delta_n}_{l_n},\Gamma_1,\cdots,
\Gamma_k,\Gamma^{\prime})\\
&&\times
\frac{(u_i)^{l_i-1}}{(l_i-1)!}[\prod_{\makebox {j=1,} j\ne i}^n\frac{(u_j)^{l_j}}{l_j!}]
[\prod_{j=1}^{k}p^{(j)}_{\Gamma_j}]W(\Delta_i,z)p_{\Gamma^{\prime}},
\end{eqnarray*}
where $\check{l_i}$ means that we omit $l_i$, and  $2h^i_--2$ is also determined by the Hurwitz formula:
\begin{equation}\label{g2}
(2-2g)d-(2-2h^i_-)=\sum_{j=1}^{n}l_j(d-l(\Delta_j))-(d-l(\Delta_i))+\sum_{j=1}^{k}(d-l(\Gamma_j))+(d-l(\Gamma^{\prime})),
\end{equation}

Moreover we have the following facts:
\begin{itemize}
\item Fact (1): By lemma 2.1,
$$W(\Delta_i,z)p_{\Gamma^{\prime}}=\frac{d!}{|C_{\Gamma^{\prime}}|}\sum_{\Gamma}z^{d+l(\Gamma^{\prime})-l(\Delta_i)-l(\Gamma)}\mu_{0}^{h^i_+,d}(\Gamma^{\prime},\Delta_i,\Gamma)p_{\Gamma};$$
\item  Fact (2):
 $\mu_{0}^{h^i_+,d}(\Gamma^{\prime},\Delta_i,\Gamma)\ne 0$ only if
\begin{equation}\label{g3}
2d-(2-2h^i_+)=(d-l(\Gamma^{\prime}))+(d-l(\Delta_i))+(d-l(\Gamma));
\end{equation}
\item  Fact (3): by the formula (\ref{g1}), (\ref{g2}), (\ref{g3}), we have
$$h=h^i_++h^i_-+l(\Gamma^{\prime})-1.$$
\end{itemize}
Then the proposition follows from Lemma 1.1.

\bigskip
\remark We note an interesting phenomenon in the above proof, i.e.,
we ``lost" the genus after we have executed the symplectic cutting:
\begin{eqnarray}
&&(2h-2)-(2h^i_--2)\nonumber\\
&=&(2h^i_+-2)+2l(\Gamma^{\prime})\nonumber\\
&=&d+l(\Gamma^{\prime}))-l(\Delta_i)-l(\Gamma),
\end{eqnarray}
which is the key observation to define the genus expanded cut-and-join operators.

 \bigskip
Immediately, we have
\corollary
\begin{equation}
\Phi_g\{z|(u_1,\Delta_1)\cdots,(u_n,\Delta_n)|p^{(1)},\cdots,p^{(k)},p\} =[\prod_{i=1}^n exp (u_iW(\Delta_i,z))]
\Phi_g\{z||p^{(1)},\cdots,p^{(k)},p\}
\end{equation}
In particular, if we take $g=0,k=0, 1$, we have
\item
$$\Phi_0\{z|(u_1,\Delta_1),\cdots,(u_n,\Delta_n)|p\}=[\prod_{i=1}^{n}exp(u_iW(\Delta_i,z))](z^{-2d}\frac{p_1^d}{d!});$$
\item
$$\Phi_0\{z|(u_1,\Delta_1),\cdots,(u_n,\Delta_n)|q,p\}
=[\prod_{i=1}^{n}exp(u_iW(\Delta_i,z))](\sum_{\Delta}z^{-2l(\Delta)}\frac{|C_\Delta|}{d!}q_\Delta p_\Delta).$$
Moreover, if we take $\Delta=(2,1,\cdots,1)\vdash d$, then
we obtain the generating functions for
the classical almost simply Hurwitz numbers and the almost simply double Hurwitz numbers
\cite{[GJ1]}, \cite{[GJ2]}, \cite{[GJ3]},\cite{[LZZ]}, \cite{[OP]}:
\begin{equation}\label{22}
\Phi_0\{z|(u,\Delta)|p\}=exp(uW((2,1,\cdots,1),z))(z^{-2d}\frac{p_1^d}{d!});
\end{equation}
\begin{equation}
\Phi_0\{z|(u,\Delta)|p,q\}=exp(uW((2,1,\cdots,1),z)(\sum_{\Delta}z^{-2l(\Delta)}\frac{|C_\Delta|}{d!}q_\Delta p_\Delta)
\end{equation}

\example
Assume $d=3$ and $\Delta=(2,1),$ then we have
$$
W(\Delta,z)=z^2[2p_3\frac{\partial^2}{\partial p_1\partial p_2}+
\frac{1}{2}p_1p_2\frac{\partial^3}{\partial p_1\partial p_1\partial p_1}]+
[p_1^3\frac{\partial^2}{\partial p_1\partial p_2}+3p_1p_2\frac{\partial}{\partial p_3}].
$$
Thus  we get the generating function(with any genus) from (\ref{22}) for $d=3$ (which also can be checked by direct
calculation):
$$
\Phi_0\{z|(u,\Delta)|p\}=\frac{1}{6}p_1^3z^{-6}+ \frac{1}{2}up_1p_2z^{-4}+\frac{1}{2}u^2p_3z^{-2}+\frac{1}{4}u^2p_1^3z^{-4}+\frac{3}{4}u^3p_1p_2z^{-2}+\cdots,
$$
where we omit the higher terms of $u.$

\bigskip
{\bf Acknowledgements} The author is thankful to his advisors: Prof. An-Min Li and Prof. Guosong Zhao,
who introduced him to the field of geometry. He would like to thank the every members of geometry
team of Mathematics College in Sichuan University, and Prof. Yongbin Ruan, Prof. Yuping Tu,
Prof. Guohui Zhao, Prof. Qi Zhang for their helpful discussions. He  is very appreciative  of  the referees' suggestions .

\newpage
\mbox{}\\[10pt]

\end{document}